\documentclass[12pt]{amsart} 
\usepackage{amssymb}

\def\marginpar#1{}
\makeatletter \@mparswitchfalse \makeatother

\newtheorem{question}{Question}
\newtheorem{theorem}{Theorem}
\newtheorem{corollary}{Corollary}

\begin{document}

% bring back \eqalign from plain TeX
\def\eqalign#1{\null\,\vcenter{\openup\jot
  \ialign{\strut\hfil$\displaystyle{##}$&$\displaystyle{{}##}$\hfil
      \crcr#1\crcr}}\,}

\let\\\cr
\let\phi\varphi
\let\union\bigcup
\let\inter\bigcap
\def\supp{\operatorname{supp}}
\let\emptyset\varnothing

\title[Invariance of the minimal envelope map of Douglas algebras]
{Some conditions on Douglas algebras that imply the invariance of the
minimal envelope map}
\author{Carroll Guillory}
\address{Department of Mathematics\\
University of Southwestern Louisiana\\
Lafayette, LA 70504}

\curraddr{Mathematical Science Research Institute\\
100 Centennial Drive\\
Berkeley, CA 94720}

\email{cjg2476@@usl.edu}

\thanks{All work on this paper was done while the author was at the
Mathematical Science Research Institute. The author thanks the
Institute for its support during this period.  Research at MSRI is
supported in part by NSF grant DMS-9022140.}

\begin{abstract}        
We give several conditions on certain families of Douglas
algebras that imply that the minimal envelope of the given algebra is the
algebra itself. We also prove that the minimal envelope of the
intersection of two Douglas algebras is the intersection of their
minimal envelope.
\end{abstract}        

\maketitle

\section{Introduction}

Let $D$ denote the open unit disk in the complex plane, and $T$ the unit
circle. By $L^\infty$ we mean the space of essentially bounded measurable
functions on $T$ with respect to the normalized Lebesgue measure. We
denote by $H^\infty$ the space of all bounded analytic functions in $D$. Via
identification with boundary functions, $H^\infty$ can be considered as a
uniformly closed subalgebra of $L^\infty$. Any uniformly closed subalgebra $B$
strictly between $H^\infty$ and $L^\infty$  is called a Douglas 
algebra. We denote by $M(B)$  
the maximal ideal space of a Douglas algebra $B$.  If $C$ is the set of
all continuous functions on $T$, we set $H^\infty+C=\{h+g: h \in 
H^\infty,\, g\in C\}$.   Then $H^\infty + C$ becomes the smallest
Douglas algebra containing $H^\infty$ properly.

The function
$$
q(z) = \prod_{n=1}^\infty \frac{|z_n|}{z_n} \frac{z-z_n}{1-\bar z_nz}
$$
is called a Blaschke product if $\sum_{n=1}^\infty (1-|z_n|)$ 
converges. The set $\{z_n\}_n$ is called the zero set of 
$q$ in $D$. Here ${|z_n|}/{z_n} = 1$ is understood whenever $z_n=0$.
We call $q$ an interpolating Blaschke product if
$$
\inf_n \prod_{m:m \neq n} \left | \frac{z_m - z_n}{1-\bar z_nz_m}
\right | > 0.
$$
An interpolating Blaschke product $q$ is called sparse (or thin) if 
$$
\lim_{n\to \infty} \prod_{m:m\neq n} \left | \frac{z_m-z_n}{1-\bar
z_nz_m} \right | = 1.
$$

The set $Z(q) = \{x \in M (H^\infty)\setminus D: q(x) = 0\}$ is called
the zero set of $q$ in $M(H^\infty + C)$. Any function $h$ in
$H^\infty$ with $|h|=1$ almost everywhere on $T$ is called an inner
function. Since $|q|=1$ for any Blaschke product, Blaschke products
are inner functions. Let $QC = (H^\infty + C)\cap 
\overline{(H^\infty + C)}
\text{ and for } x \in M (H^\infty + C)$, set 
$$
Q_x = \{y \in M (L^\infty): f(x) = f(y) \text{ for all } f\in 
QC\}.
$$
Then $Q_x$ is called the $QC$-level set for $x$. For $x\in M(H^\infty +
C)$, we denote by $\mu_x$ the representing measure for $x$, and its
support set by $\supp\mu_x$. By $H^\infty[\bar q]$ we mean the
Douglas algebra generated by $H^\infty$ and the complex conjugate of
the function $q$. Since $M(L^\infty)$ is the Shilov boundary for every
Douglas algebra, a closed set $E$ contained in $M(L^\infty)$ is called
a peak set for a Douglas algebra $B$ if there is a function $f$ in  $B$ with
$f=1$ on $E$ and $|f| < 1$ on $M(L^\infty)\setminus E$. A closed set
$E$ is a weak peak set for $B$ if $E$ is the intersection of a family
of peak sets. If the set $E$ is a weak peak set for $H^\infty$ and we
define
$$
H_E^\infty =\{f\in L^\infty: f|_E \in H^\infty|_E\},
$$
then $H_E^\infty$ is a Douglas algebra. For a Douglas algebra $B$,
we define $B_E$ similarly.

For an interpolating Blaschke product $q$ we put $N(\bar q)$ the
closure of 
$$
\bigcup\{\supp\mu_x:x\in M(H^\infty + C),\, |q(x)| <1\}.
$$ 
Then
$N(\bar q)$ is a weak peak set for $H^\infty$. By $N_0(\bar q)$
we denote the closure of $\bigcup \{\supp\mu _x:x\in Z(q)\}$. In
general $N_0(\bar q)$ does not equal $N(\bar q)$, but in this paper
$N_0(\bar q) = N(\bar q)$. For $x\in M(H^\infty)$, we let $E_x = \{y\in
M(H^\infty): \supp\mu_y = \supp\mu_x\}$ and call $E_x$ the level set
of $x$. Since the sets $\supp\mu_x$ and $N(\bar q)$ are weak
peak sets for $H^\infty$, both $H^\infty_{\supp\mu_x}$ and
$H^\infty_{N(\bar q)}$ are Douglas algebras. For the interpolating
Blaschke product $q$, set
$$
A=\bigcap_{x\in M(H^\infty + C)} \{H^\infty_{\supp\mu_x}: |q(x)| < 1\}
$$
and
$$
A_0 = \bigcap\{H^\infty_{\supp\mu_x}:x \in Z(q)\}.
$$
Our assumptions on $q$ throughout this paper imply that
$H^\infty_{N(\bar q)} = A = A_0$ (see \cite{10}). For $x$ and $y$ in
$M(H^\infty)$, the pseudohyperbolic distance $\rho$ is defined by
$$
\rho(x,y)=\supp\{|h(y)|:\|h\|_\infty\leq 1,\, h \in H^\infty,\, h(x)=0\}.
$$
For $x$ and $y$ in $D$, we have 
$$
\rho(x,y)=\left | \frac{x-y}{1-\bar y x} \right |.
$$ 
For $x \in M(H^\infty)$, we define the Gleason part
$P_x$ of $x$ by 
$$
P_x = \{y \in M(H^\infty): \rho(x,y)< 1\}.
$$
If $P_x \neq\{x\}$, then $x$ is said to be a nontrivial point. 
For the definition of those interpolating Blaschke 
products that are of type G and of finite type G, see \cite{8}. 

Let $B$ be a Douglas algebra. The Bourgain algebra $B_b$ of $B$
relative to $L^\infty$ is the set of $f$ in $L^\infty$ such that
$\|ff_n +B\|_\infty \to 0$ for every sequence $\{f_n\}$ in $B$ with
$f_n \to 0$ weakly in $B$. An algebra $B$ is called a minimal
superalgebra of an algebra $A$ if $A \subset B$ and  $\supp\mu_x = \supp\mu_y$
for all $x,y \in M(A)\setminus M(B)$. The minimal envelope $B_m$ of a
Douglas algebra $B$ is defined to be the smallest Douglas algebra
that contains all the minimal superalgebras of $B$. The mapping that
assigns to $B$ the Douglas algebra $B_m$ is called the minimal
envelope map.

\section{Conditions for $B_m=B$}

We begin with the following theorem. The case for the Bourgain algebra
$B_b$ has been proven in \cite{15}. The proof used here is quite
different from theirs, and can be used to show that this result also
holds for the Bourgain algebras.

\begin{theorem}
Let $A$ and $B$ be Douglas algebras. Then $(A\cap B)_m = A_m\cap
B_m$.
\end{theorem}

\begin{proof}
Since $A\cap B$ is contained in both $A$ and $B$ by Proposition 6 of
\cite{11}, $(A\cap B)_m$ is contained in both $A_m$ and $B_m$. Hence
$(A\cap B)_m \subset A_m \cap B_m$.

To show that $A_m\cap B_m \subset (A\cap B)_m$, let $\psi$ be an
interpolating Blaschke product such that $\bar\psi \in A_m \cap B_m$.
We show that $\bar\psi \in (A\cap B)_m$. By Theorem D of \cite{11} we
can assume that there is an $x\in M(A)$ and a $y\in M(B)$ such that
$\{\lambda \in M(A):|\psi(\gamma)| < 1\}=E_x$ and
$\{\mu \in M(B):|\psi(\mu)| < 1\} = E_y$. This implies that 
$M(A)=M(A[\bar\psi])\cup E_x$ and $M(B)=M(B[\bar\psi])\cup E_y$. Thus
$$
\eqalign{
M(A\cap B) &= M(A)\cup M(B)
            = M(A[\bar\psi])\cup M(B[\bar\psi])\cup E_x \cup E_y\\
           &= M(A[\bar\psi]\cap B[\bar\psi])\cup E_x \cup E_y
            = M((A\cap B)[\bar\psi])\cup E_x \cup E_y.
}
$$
Hence $\{ \omega \in M(A\cap B):|\psi(\omega)|<1\}=E_x \cup E_y$. By
Theorem 3 of \cite{11} we have $\bar{\psi}\in(A\cap B)_m$. Our
theorem follows.
\end{proof}

The following corollaries are immediate consequences of the theorem.

\begin{corollary}
Let $A$ and $B$ be Douglas algebras with $A=A_m$ and $B=B_m$. Then
$A\cap B=(A\cap B)_m$.
\end{corollary}

\begin{corollary}
Let $B_0$ be a Douglas algebra with $(B_0)_m=B_0$ and $B$ be a Douglas
algebra such that there is an interpolating Blaschke product $q$ with
$\{\lambda \in M(B):|q(\lambda)|=1\}\subset M(B_0)$. If $A$ has the property
that $A_m=A$, then $B[\bar q]\cap A = (B[\bar q]\cap A)_m$. In
particular, $H^\infty[\bar\phi]\cap A = (H^\infty[\bar q]\cap A)_m$ for
every interpolating Blaschke product $\phi$.
\end{corollary}

\begin{proof}
The hypothesis $\{\lambda \in M (B): |q(\lambda)|=1\}\subset M(B_0)$
implies that $B_0[\bar q] = B[\bar q]$, hence $[B[\bar q]]_m = B[\bar
q]$. So $B[\bar q]\cap A = (B[\bar q]\cap A)_m$ follows from our
theorem. By Theorem 1 of \cite{2} we have $H^\infty_m = H^\infty + C$, so
$\{\lambda \in M(H^\infty): |\psi(\lambda)|=1\} \subset M (H^\infty +
C)$ for every interpolating Blaschke product. The second part of the
corollary follows.
\end{proof}

\begin{theorem}
Let $B_0$ be a Douglas algebra such that $(B_0)_m = B_0$. Let $B$ be
any other Douglas algebra such that there is an interpolating Blaschke
product $\phi$ with $\{\lambda \in M(B):|\phi(\lambda)|=1\} \subset
M(B_0)$. For an interpolating Blaschke product $q$ set $B_x =
B[\bar\phi]\cap H^\infty_{\supp\mu_x}$ for each $x \in Z(q)$. Put $B_e
= \bigcap\{B_x: x \in Z(q)\}$.

\begin{enumerate}
\item[(i)]  If $q$ is of finite type $G$, then $(B_e)_b = B_e$.
\item[(ii)] If $q$ is the product of a finite number of sparse interpolating
Blaschke products, then $(B_e)_m = B_e$. 
\end{enumerate}
\end{theorem}

\begin{proof}
We show that $B_e = B[\bar\phi]\cap H^\infty_{N(\bar q)}$. By,
an unpublished result of D. Sarason,
$$
\eqalign{
M(B_e) &= M\Bigl(\bigcap_{x \in Z(q)}B_x\Bigr)
= \overline{\bigcup_{x\in Z(q)}M(B_x)}\\ &
= \overline{\bigcup_{x\in Z(q)}M(B[\bar\phi]\cap H^\infty_{\supp\mu_x})}  
=\overline{\bigcup_{x\in Z(q)}(M(B[\bar\phi])\cup
M(H^\infty_{\supp\mu_x}))}\,,\\
M(B_e)&=M(B[\bar\phi])\cup\overline{\bigcup_{x\in
Z(q)}M(H^\infty_{\supp\mu_x})}\,.}
$$

Now, if $q$ is of finite type $g$, Proposition 1 of \cite{9}
and Theorem 3.2(i) of \cite{8} imply
$
\bigcap_{x\in Z(q)} H^\infty_{\supp\mu_x} = H^\infty_{N(\bar q)}.
$ 
Thus
$$
\overline{\bigcup_{x\in Z(q)}M(H^\infty_{\supp\mu_x})}=M(H^\infty_{N(\bar q)})
$$ 
and we get
$$
\eqalign{
M(B_e) &= M(B[\bar\phi])\cup M(H^\infty_{N(\bar q)})\\
       &= M(B[\bar\phi]\cap H^\infty_{N(\bar q)}).
}
$$
So by the Chang--Marshall Theorem \cite{1, 13} we have $B_e =
B[\bar\phi]\cap H^\infty_{\supp\mu_x}.$

The condition $\{\lambda\in M(B): |\phi(\lambda)|=1\} \subset M( B_0)$
implies that $B[\bar\phi]=B_0[\bar\phi]$. Hence $[B[\bar\phi]]_m =
[B_0[\bar\phi]]_m = [B_0]_m[\bar\phi]=B_0[\bar\phi]=B[\bar\phi]$,
where the middle equality follows from Theorem 4 of \cite{11}. By
theorem 1 of \cite{7} (and its proof) we have $(H^\infty_{N(\bar
q)})_b=H^\infty_{N(\bar q)}$ if $q$ is of finite type $G$. Thus
$$
(B_e)_b=(B[\bar\phi]\cap H^\infty_{N(\bar q)})_b =
B[\bar\phi]\cap H^\infty_{N(\bar q)} = B_e.
$$ 
This proves (i).

Now if $q$ is the product of a finite number of sparse Blaschke
products, then again by Theorem 1 of \cite{7} we have $(H^\infty_{N(\bar
q)})_m=H^\infty_{N(\bar q)}$, and so $(B_e)_m=B_e$. This proves
(ii).
\end{proof}
 
\begin{corollary}
Suppose the hypothesis of Theorem 2 holds if $\phi=q$. Then $B_e=B$, and:
\begin{enumerate}
\item[(i)] $B = (B)_b$ if $q$ is of finite type $G$; 
\item[(ii)] 
$B = B_m$ is the product of a finite number of sparse Blaschke products.
\end{enumerate}
\end{corollary}

\begin{theorem}

Let $B$ be a Douglas algebra and let $E$ be a peak set for $B$. Then
$(B_E)_m = B_E$.

\end{theorem}
\begin{proof}
Since $E$ is a peak set for $B$, there is a function $f \in B$ such
that $f=1$ on $E$ and $|f|< 1$ on $M(L^\infty)\setminus E$. By
\cite[p.~39]{3} we have 
$$
M(B_E) = \{x\in M(B):\supp\mu_x \subseteq E\}\cup M(L^\infty)
       = \{x\in M(B): f(x)=1\}\cup M(L^\infty).
$$
Let $I$ be any interpolating Blaschke product such that $Z(I)\cap
M(B_E)\neq\emptyset$.  We will show that there is an uncountable set
$\Gamma\subset Z(I)\cap M(B_E)$ such that $E_\alpha\neq E_\beta$ for
distinct $\alpha,\beta\in\Gamma$. Let $\{z_n\}$ be the zero sequence
of $I$ in $D$ and $y\in Z(I)\cap M(B_E)$. Since $y\in Z(I)$, there is
a subsequence $\{z_{n_k}\}$ of $\{z_n\}$ such that $z_{n_k}\to y$.
Since $y\in M(B_E)$, we have $f(z_{n_k})\to 1$.  Let $I_1$ be the
factor of $I$ such that $\{z_{n_k}\}$ is the zero sequence of $I_1$.
Then
$$
Z(I_1) = \overline{\{z_{n_k}\}}\setminus\{z_{n_k}\}
$$
and since $Z(I_1)$ is equivalent to the \v Cech compactification
of the integers, $Z(I_1)$ is an uncountable set. Since $f(z_{n_k})\to
1$, we have that $f=1$ on $Z(I_1)$. Thus $Z(I_1)\subset Z(I)\cap
M(B_E)$, and so $Z(I)\cap M(B_E)$ is an infinite set. To show that
$\Gamma$ exists, take a subsequence
$\{z_{n_{k_0}}\}$ of $\{z_{n_k}\}$ such that $\{z_{n_{k_0}}\}$ is a
sparse Blaschke sequence in $D$ (see \cite{4}). Let $I_0$ be the
sparse Blaschke product whose zero sequence in $D$ is
$\{z_{n_{k_0}}\}$. Then, by Theorem 1 of \cite{12}, we have 
$$
N(I_0)=\bigcup_{x\in Z(I_0)}Q_x.
$$
Since $Z(I_0)$ is uncountable, Lemma 4 of \cite{12} says that
$Q_x\neq Q_y$ for distinct $x,y \in Z(I_0)$. This
implies that $\supp\mu_x\cap\supp\mu_y=\emptyset$ and $E_x\neq E_y$.
Set $\Gamma = Z(I_0)$. Then 
$$
\bigcup_{\alpha\in\Gamma}E_\alpha\subseteq\{\lambda\in M(B_E):|I(\lambda)| <
1\}.
$$
By Theorem 3 of \cite{11}, $\bar{I}\not\in(B_E)_m$. This shows
that $(B_E)_m = B_E$.
\end{proof}

\begin{theorem}

Let $A$ be a Douglas algebra and let $\{q_n\}$ be a sequence of
interpolating Blaschke products such that $\overline{q_n}\not\in A$
for every $n$. Let $B = A[\overline{q_n}; n=1,2,3\ldots]$. Then
$B_m=B$. 
\end{theorem}

\begin{proof}
Since $M(B)=\{\lambda\in M(A):|q_{n(\lambda)}|=1, \,n=1,2,\ldots\}$, we
can define the function $F$ on $M(H^\infty)$ as
$$
F(x)=\sum^\infty_{n=1}\frac{1}{2^n}|q_n(x)|.
$$
Then $F=1$ on $M(B)$ and $|F|<1$ on $M(H^\infty)\setminus M(B)$. Let
$C$ be any interpolating Blaschke product such that $Z(C)\cap
M(B)\neq\emptyset$. Let $\{z_n\}$ be the zero sequence of $C$ in $D$.
As in the proof of Theorem 3 we can find a sparse subsequence
$\{z_{n_k}\}$ of $\{z_n\}$ such that $F=1$ on
$\overline{\{z_{n_k}}\}\setminus\{z_{n_k}\}$. If we set $C_1$ to be
the factor of $C$ whose zero sequence in $D$ is $\{z_{n_k}\}$, then as
in the proof of Theorem 3 we can show that $E_{\alpha}\neq
E_{\beta}$ for distinct $\alpha, \beta\in Z(C_1)$, and
$$
\bigcup_{\alpha\in Z(C_1)}E_{\alpha}\subset\{\lambda\in  
M(B):|C \lambda | < 1\}. 
$$
So, by Theorem 3 of \cite{11}, we have 
$\bar C\not\in B_m$. This implies that $B_m = B$. 
\end{proof}

\begin{theorem}
Let $\{q_n\}$ be a sequence of interpolating Blaschke products not
invertible in a Douglas algebra $A$. Suppose that for each $n$ there is a
Douglas algebra $A_n$ with $(A_n)_m = A_n$ and
$$
\{\lambda\in M(A):|q_n(\lambda)|=1\}\subset M(A_n).
$$
Let $B=\bigcap^\infty_{n=1}A[\bar q_n]$. Then $B_m=B$.
\end{theorem}

\begin{proof}
Since $B\subseteq A[\bar q_n]$ for $n=1, 2 \ldots$, by
Proposition 6 of \cite {11} we have
$B_m\subseteq(A[\bar q_n])_m$. Hence
$$
B_m\subseteq\bigcap^\infty_{n=1}(A[\overline {q}_n])_m.
$$ 
To show that
$B_m=B$ it suffices to show that
$(A[\bar q_n])_m=A[\bar q_n]$ for all $n$. Since
$\{\lambda\in M(A): |q_n(\lambda)|=1\}\subset M(A_n)$, we have
$A[\bar q_n]=A_n[\bar q_n]$. So using Theorem 4 of
\cite{11} we have
$$
(A[\bar q_n])_m = (A_n[\bar q_n])_m
                      = (A_n)_m[\bar q_n]
                      = A_n[\bar q_n]
                      = A[\bar q_n].
$$
Thus $B_m\subseteq\bigcap^\infty_{n=1}(A[\bar q_n])_m
=\bigcap^\infty_{n=1}A[\bar q_n]=B$.
\end{proof}

\begin{corollary}
If $B=\bigcap^\infty_{n=1}H^\infty[\bar q_n]$, then $B_m=B$.
\end{corollary}

\begin{theorem}
Let $B$ be a Douglas algebra. Then $B\subset B_m$ if and only if there
exists a point $x\in M(B)$ whose level set $E_x$ is an open
subset in $M(B)$.
\end{theorem}

\begin{proof}
Let $B\subset B_m$. By Theorem D of \cite {11} there exists an
interpolating Blaschke product $\phi$ such that $\{\lambda\in
M(B):|\phi(\lambda)| < 1 \} = E_x$ for some $x\in M(B)$. By Lemma
9(ii) of \cite{11}, $E_x$ is an open subset of $M(B)$.

Suppose that $B=B_m$. Then for every $y\in M(B)$ and every
interpolating Blaschke product $q$ with $q(y)=0$ the set $Z(q)\cap
M(B)$ contains an infinite subset $\Gamma$ such that $E_\alpha\neq
E_\beta$ for all distinct $\alpha, \beta\in\Gamma$, and
$$
\bigcup_{\alpha\in\Gamma} E_\alpha = \{\lambda\in M(B):|q(\lambda)| < 1\}.
$$
We'll show that there is a subset $\{x_\sigma\}_{\sigma\in A}$ of
$\Gamma$ such that 
$y\in\overline{\{x_\sigma\}}_{\sigma\in A}\setminus\{y\}$, 
that is,
$$
E_y\subseteq\overline{\{E_\alpha\}}_{\alpha\in\{\Gamma\setminus\{y\}\}}.
$$
Suppose to the contrary that
$E_y\not\subset\overline{\{E_\alpha\}}_{\alpha\in\{\Gamma\setminus\{y\}\}}$. 
Then $\overline{\{E_\alpha\}}_{\alpha\in\{\Gamma\setminus\{y\}\}}$
and $\{y\}$ are two closed disjoint subsets of $M(H^\infty)$. Choose 
neighborhoods $V$ of 
$\overline{\{E_\alpha\}}_{\alpha\in\{\Gamma\setminus\{y\}\}}$
and $U$ of $\{y\}$ such that $\bar U\cap\bar V=\emptyset$.
Let $q_1$ be the the factor of $q$ with zeros in $U\cap D$.
Then we see that $Z(q_1)\cap M(B)\subset E_y$ and $\{\lambda\in
M(B):|q_1(\lambda)|< 1\} = E_y$. Thus $\bar{q}_1\in B_m\setminus 
B$, which is a contradiction. Hence 
$$
E_y\subset\overline{\{E_\alpha\}}_{\alpha\in\{\Gamma\setminus\{y\}\}}.
$$
Thus there is a subset $\{x_\sigma\}_{\sigma\in A}\subset\Gamma$ such
that $x_\sigma \to y$, and $x_\sigma\not\in E_y$ for all $\sigma\in
A$. This implies that $M(B)\setminus E_y$ is not a closed subset of
$M(H^\infty)$. Thus $E_y$ is not an open subset of $M(B)$. This proves
our theorem.
\end{proof}

I have been unable to answer the following two questions, which I
close with.

\begin{question} 
Find an $x\in M(H^\infty + C)\setminus M(L^\infty)$ such that
$H^\infty_{\supp\mu_x} = (H^\infty_{\supp\mu_x})_m$ (if such an $x$
exists).
\end{question}

\begin{question}
Find a Douglas algebra $B$ such that $B\subsetneq B_b\subsetneq B_m$
(if such a Douglas algebra exists).
\end{question}

An immediate consequence of Theorem 6 above and Theorem 7 of \cite{14}, 
related to Question 2, is the following:

\begin{theorem}
Let $B$ be a Douglas algebra. Then $B=B_b\subsetneq B_m$ if and only
if for all $x\in M(B)$ the set $P_x$ is not open in $M(B)$ but for
some $x\in M(B)$ the set $E_x$ is an open subset in $M(B)$.
\end{theorem}

\end{document}